\documentstyle[twoside,psfig,amssymb]{article}

\oddsidemargin=\evensidemargin
\catcode`\@=11
\long\def\@makefntext#1{
\protect\noindent \hbox to 3.2pt {\hskip-.9pt  
$^{{\eightrm\@thefnmark}}$\hfil}#1\hfill}		

\def\ps@myheadings{\let\@mkboth\@gobbletwo		
\def\@oddhead{\hbox{}
\rightmark\hfil\eightrm\thepage}   
\def\@oddfoot{}\def\@evenhead{\eightrm\thepage\hfil
\leftmark\hbox{}}\def\@evenfoot{}
\def\sectionmark##1{}\def\subsectionmark##1{}}

\catcode`@=11		      		     
\def\ps@plain{\let\@mkboth\@gobbletwo
     \def\@oddhead{}\def\@oddfoot{\eightrm\hfil\thepage
     \hfil}\def\@evenhead{}\let\@evenfoot\@oddfoot}

\newcommand{\nonumsection}[1] {\vspace{12pt}\noindent{\tenbf #1}
	\par\vspace{5pt}}

\newcommand{\textlineskip}{\baselineskip=13pt}
\newcommand{\smalllineskip}{\baselineskip=10pt}

\newcommand{\copyrightheading}[1]
	{\vspace*{-2.5cm}\smalllineskip{\flushleft
	{\footnotesize Journal of Knot Theory and Its Ramifications #1}\\
   	{\footnotesize \copyright\kern2pt World Scientific 
         Publishing Company}\\
         }}

\renewenvironment{thebibliography}[1]
	{\frenchspacing
	 \ninerm\baselineskip=11pt
	 \begin{list}{[\arabic{enumi}]}
	{\usecounter{enumi}\setlength{\parsep}{0pt}
	 \setlength{\leftmargin 13.7pt}{\rightmargin 0pt} 
	 \setlength{\itemsep}{0pt} \settowidth
	{\labelwidth}{[#1]}\sloppy}}{\end{list}}

\def\pmb#1{\setbox0=\hbox{#1}
	\kern-.025em\copy0\kern-\wd0
	\kern.05em\copy0\kern-\wd0
	\kern-.025em\raise.0433em\box0}

\def\fpage#1{\begingroup
\voffset=.3in
\thispagestyle{empty}\begin{table}[b]\centerline{\footnotesize #1}
	\end{table}\endgroup}

\def\runninghead#1#2{\pagestyle{myheadings}
\markboth{{\protect\footnotesize\it{\quad #1}}\hfill}
{\hfill{\protect\footnotesize\it{#2\quad}}}}

\def\abstracts#1#2#3#4{{
	\centering{\begin{minipage}{4.5in}\footnotesize\baselineskip=10pt
	\centerline{ABSTRACT} 
	\parindent=15pt #1\par 
	\parindent=15pt #2\par
	\parindent=15pt #3\par
	\parindent=15pt #4\par
	\end{minipage}}\par}} 

\font\tenbf=cmbx10

\font\ninerm=cmr9
\font\nineit=cmti9
\font\ninebf=cmbx9
\font\eightrm=cmr8

\def\qed{\hbox{${\vcenter{\vbox{			
   \hrule height 0.4pt\hbox{\vrule width 0.4pt height 6pt
   \kern5pt\vrule width 0.4pt}\hrule height 0.4pt}}}$}}

\begin{document}
\setlength{\textheight}{7.7truein}  

\runninghead{Unusual formulae for the Euler characteristic}
{Unusual formulae for the Euler characteristic} 

\normalsize\textlineskip
\thispagestyle{empty}
\setcounter{page}{1}

\copyrightheading{}		    

\vspace*{0.88truein}

\fpage{1}
\centerline{\bf UNUSUAL FORMULAE FOR THE EULER CHARACTERISTIC}
\baselineskip=13pt
\vspace*{0.37truein}
\centerline{\footnotesize JUSTIN ROBERTS}
\baselineskip=12pt
\centerline{\footnotesize\it Department of Mathematics}
\baselineskip=10pt
\centerline{\footnotesize\it UC San Diego}
\centerline{\footnotesize\it 9500 Gilman Drive}
\centerline{\footnotesize\it La Jolla CA 92093}
\centerline{\footnotesize\it USA}

\vspace*{1in}

\abstracts{\noindent Everyone knows that the Euler characteristic of a combinatorial 
manifold is given by the alternating sum of its numbers of
simplices. It is shown that there are other linear combinations of the
numbers of simplices which are combinatorial invariants, but that all
such invariants are multiples of the Euler characteristic.}{}{}{}

\textlineskip			
\vspace*{12pt}			


\noindent{\bf Theorem [Pachner 1987, \cite{[3]}].}
Two closed combinatorial $n$-manifolds are $PL$-homeomorphic if and
only if it is possible to move between their triangulations using a
sequence of {\em Pachner moves} {\rm (}bistellar moves{\rm )} and
simplicial isomorphisms.

\vspace*{5pt}

The Pachner moves consist of replacing one part of the boundary of an
$(n+1)$-simplex by the other; specifically, one views the simplex as a
join of simplices $A * B$, and replaces $\partial A * B$ with $A *
\partial B$. Thus in an $n$-manifold there are only $\lfloor
\frac{n}{2}\rfloor +1$ different kinds of moves needed.

The natural inclination of the topologist on seeing this theorem is to
look for quantities invariant under the Pachner moves so as to obtain
invariants of combinatorial manifolds. The Turaev-Viro invariant
\cite{[4]} of $3$-manifolds can be proved to be invariant using exactly
this method, but the most obvious example is clearly the Euler
characteristic, defined as the alternating sum of numbers of
simplices. One can simply write down how the numbers of simplices
alter under the Pachner moves and check invariance directly.

It is obvious that there are many {\em other} linear combinations of
the functions $f_0, f_1, \ldots, f_n$ ($f_i$ being the number of
$i$-simplices) which are combinatorial invariants, purely because
invariance under the Pachner moves imposes only about $\frac{n}{2}$
conditions. What are these invariants? 

Now in fact Pachner goes some way towards answering this question in
his original paper. He observes that if one defines the $h$-vector
$(h_0, h_1, \ldots, h_{n+1})$ of a combinatorial manifold in the same
way as one does for convex polytopes, where each $h_i$ is a certain
linear combination of the $f_j$'s and a constant term (see Ziegler
\cite{[5]} for the precise formulae), then the $\lfloor \frac{n+1}{2}
\rfloor$ quantities $h_i - h_{n+1-i}$ turn out to be invariants under
the Pachner moves. Consequently one can write an explicit set of
linear combinations of the $f_i$ which in fact spans the space of
invariants ($\lfloor \frac{n}{2}\rfloor +1 + \lfloor \frac{n+1}{2}
\rfloor = n+1$).

However, he gives no hint as to what these invariants actually are,
though it should be clear to any topologist that they will be
expressible in terms of classically known invariants, by the ``no free
lunch'' principle. (On the evidence of how many there seem to be, one
might hope that they are the Betti numbers.) The only evaluation he
mentions is that for any combinatorial $n$-sphere, the $h_i-h_{n+1-i}$
are all zero: these identities are the classical {\em Dehn-Sommerville
equations}, giving constraints on the numbers of simplices in the
boundary of a convex polytope.

Let us now think about the properties of such invariants. Any linear
combination $v$ of the $f_i$ which is an invariant of combinatorial
$n$-manifolds will obviously be multiplicative under finite coverings,
and additive in the same sense as the Euler characteristic (the
``valuation'' property $v(A \cup B) = v(A) + v(B) - v(A \cap
B)$). This is why the invariants can't be Betti numbers, and why it's
difficult to imagine that they could be anything other than the Euler
characteristic, as in fact turns out to be the case.

\vspace*{5pt}

\noindent{\bf Theorem.}
Any linear combination $v$ of the numbers of simplices which is an
invariant of closed combinatorial $n$-manifolds is proportional to the
Euler characteristic. {\rm (}Thus in odd dimensions it is identically zero,
and in even dimensions the constant of proportionality is
$\frac{1}{2}v(S^n)$.{\rm )}

\vspace*{5pt}

After finding a nice short proof of this theorem I then trawled the
literature, certain that it must be known. In fact I could not find
any explicit description of this fact. The literature on polytopes
seems to concentrate on the Dehn-Sommerville equations for
combinatorial spheres, without addressing arbitrary combinatorial
manifolds. Eventually I discovered MacDonald \cite{[2]}, who proves a
theorem amounting to the Dehn-Sommerville equations on a general
combinatorial manifold which indeed implies the above result.  The
topological part of the proof is very short, but to get to the above
result one then has to use the fact that the $h_i-h_{n+1-i}$
invariants span the space under consideration, which requires a bit of
work with binomial coefficients. So here is an alternative.

\vspace*{5pt}

\noindent{\bf Proof.}
Assume we have a linear combination $v$ of the $f_i$ which is an
invariant for closed combinatorial $n$-manifolds. It will not be
invariant for $n$-manifolds with boundary, but the modified version
$\tilde v(M) = v(M) -\frac{1}{2}v(\partial M)=\frac12 v(DM)$, where
$DM$ is the double of $M$, clearly will be.

Given a closed combinatorial $n$-manifold $M$, let us compute $\tilde
v(M) = v(M)$ by taking a combinatorial handle decomposition of $M$ and
examining how $\tilde v$ changes as we add handles one at a
time. (Such handle decompositions always exist: take regular
neighbourhoods of the barycentres of the simplices of $M$. See
Lickorish \cite{[1]} for an extensive discussion of this kind of
handlebody theory.)

Each $0$-handle contributes $\tilde v(B^n)$, whatever its
triangulation. Suppose now that we attach an $i$-handle $H$ to an
existing manifold-with-boundary $N$, and that the boundary of $H$ is
decomposed as $\partial_-H \cup_C
\partial_+H$, with $\partial_-H$ the part glued to $N$ and $C$
standing for ``corner''.  It is obvious that 
\[v(\partial H) = v(\partial_-H) + v(\partial_+H) - v(C)\]
and thus that
\begin{eqnarray*}  
\tilde v(N \cup H) & = &v(N \cup H) - \frac12 v(\partial (N \cup H))\\ 
&=& (v(N)+v(H)-v(\partial_-H))
-\frac{1}{2}(v(\partial N) -v(\partial_-H) +v(\partial_+H))\\ 
&= &\tilde v(N) + \tilde v(B^n) -\frac12 v(C).
\end{eqnarray*}
The same computation for the double of the handle $H$ along
$\partial_-H$ shows that
\[\tilde v (S^i \times B^{n-i}) = \tilde
v(B^n) + \tilde v(B^n) -
\frac12 v(C).  \]

Consequently the addition of an $i$-handle alters $\tilde v(N)$ by the
fixed quantity $\tilde v(S^i \times B^{n-i})-\tilde v(B^n)$,
independently of how the handle is attached, and we see that $\tilde
v(M)$ must be a linear combination of the numbers of $i$-handles
$m_i$. Because it is a combinatorial invariant, it must be invariant
under addition of a cancelling handle pair in any dimension; therefore
the coefficients of $m_i, m_{i+1}$ must be equal and opposite for all
$i$, so that $v$ is proportional to the Euler characteristic. \qed\kern0.8pt 

\vspace*{5pt}

One can interpret this as giving universal relations between the
numbers of simplices in combinatorial manifolds.  In an $n$-manifold
one always has the relation $(n+1)f_n - 2 f_{n-1}=0$, and the above
theorem says essentially that there are lots of other linear
constraints of this type.

It is easy to check that for $n$ odd, the Pachner moves actually do
impose $\lfloor \frac{n}{2} \rfloor$ independent conditions, so that
there are $\lfloor \frac{n}{2} \rfloor +1$ independent invariant
linear combinations of the $f_i$, including the Euler
characteristic. These all evaluate to zero, so can be treated as
constraints.

In even dimensions, the Pachner moves impose $\frac{n}{2}$ independent
conditions, so that there is an $\frac{n}{2}$-dimensional space of
invariant linear combinations. Each is a multiple of the Euler
characteristic, and the constant of proportionality defines a
surjective linear map to ${\mathbb R}$, whose kernel is an
$(\frac{n}{2}-1)$-dimensional space of universal constraints.

\nonumsection{Acknowledgements}
Partial support provided by NSF grant number DMS-0103922.

\nonumsection{References}

\end{document}